\documentclass{article}
\usepackage[T2C]{fontenc}
\usepackage[tbtags]{amsmath}
\usepackage{amsfonts,amssymb,mathrsfs,amscd}
\newtheorem{theorem}{Theorem}
\newtheorem{lemma}{Lemma}
\newtheorem{problem}{Problem}

\allowdisplaybreaks

\begin{document}

\title{Groups generated by involutions, numberings of posets, and central measures}

\author{A.~M.~Vershik\thanks{%
St.~Petersburg Department of Steklov Institute of Mathematics, St.~Petersburg State University, Institute for Information Transmission Problems.
E-mail: vershik@pdmi.ras.ru. Supported by the RSF grant 21-11-00152.}}

\date{01.07.2021}

\maketitle


\subsection{Definitions}
\label{ssec1}

An infinite countable ordered set
$\{P,\succ,\varnothing\}$ with minimal element~$\varnothing$ and no maximal elements is called a locally finite poset if all its principal ideals are finite. A monotone numbering of~$P$ (or of a part of~$P$) is an injective map  ${\phi\colon{\mathbb N}\to P}$ from the set of positive integers to~$P$ satisfying the following conditions: if $\phi(n)\succ \phi(m)$, then $n>m$;
$\phi(0)=\varnothing$. The distributive lattice~$\Gamma_P$ of all finite ideals of a locally finite poset
$\{P,\succ\}$ forms an
$\mathbb N$-graded graph (the Hasse diagram of the lattice). A monotone numbering of~$P$ is identified in a natural way with a maximal path in the lattice~$\Gamma_P$. The set~$T_P$ of all monotone numberings of~$P$, i.e., the space of infinite paths in the graph~$\Gamma_P$, can be endowed with a natural structure of a Borel and topological space. In the terminology related to the Young graph, the poset~$P$ is the set of $\mathbb{Z}_+^2$-finite ideals, i.e., Young diagrams, and monotone numberings are Young tableaux.

Let $P$ be a finite ($|P|<n\in \mathbb N$) or locally finite
($n=\infty$) poset; for each $i<n$, we define an involution~$\sigma_i$ on the space of numberings
$T_P=\{\phi\}$ of~$P$:
\begin{gather*}
(\sigma_i(\phi))(m)=\phi(m),\qquad m\ne i,i+1;
\\
(\sigma_i(\phi))(i)=\phi(i+1),\qquad
(\sigma_i(\phi))(i+1)=\phi(i)  \quad\text{if\ \ } \phi(i)\nprec
\phi(i+1).
\end{gather*}

The finite or countable group $G_P=\langle\sigma_i,\ i=1,2,\dots\rangle$ generated by this family of involutions will be called the
\textit{group of symmetries of numberings of the poset~$P$}, or the group of automorphisms of the path space of the graph~$\Gamma_P$. The finite (or infinite) symmetric group~$S_n$
(respectively,~$S_{\mathbb N}$) arises in our definition for the Young diagram~$(n-1,1)$ (respectively, $({\mathbb N},1)$).

This new class of finite and locally finite groups is of interest both from topological and dynamical points of view. It is meaningful even in the case where $P$~is an infinite diagram with two rows.

For a locally finite poset~$P$, the group~$G_P$ is a locally finite group, i.e., the inductive limit of finite groups of the same class. These groups are generalizations of the infinite symmetric group, but, in general, are not isomorphic to it.

The algebraic description of the groups~$G_P$ is as follows.

\begin{lemma}\label{lem1}
For every locally finite poset, the relations between the generators~$\sigma_i$,
$i=1,2,\dots$, defined above are as follows:
{\rm(i)} $\sigma_i^2=\{\sigma_i\cdot
\sigma_j\}^2=\operatorname{Id}$, $|i-j|>1$, $i,j=1,\dots,n-1$;
{\rm(ii)} for every $i=1,\dots,n-1$, the group generated by the elements~$\sigma_i$,~$\sigma_{i+1}$ \textup(in the classical case, it is isomorphic to~$S_3$: $\{\sigma_i \cdot
\sigma_{i+1}\}^3=\operatorname{Id}$\textup) is isomorphic to the sum of a certain \textup(possibly zero\textup) number of groups
$S_3$,~$\mathbb  Z_2$,~$\mathbb Z_3$, and their products; in particular,
$\{\sigma_i \cdot
\sigma_{i+1}\}^6=\operatorname{Id}$.
\end{lemma}

It is not clear whether the conditions imposed on the group by the lemma are sufficient for it to be the  group of symmetries of numberings of some poset~$P$. Apparently, the structure and properties of the group~$G_P$
for  a poset~$P$ have not been studied; it is of much interest to describe all such finite and countable groups up to isomorphism. For some examples, see~\cite{1}.

\subsection{Central measures on the space of numberings}
\label{ssec2}
A Borel probability measure on the set
$T_P=\{\phi\}$ of numberings of a poset~$P$ is said to be central if it is invariant under the group~$G_P$. Every central measure determines a random numbering of the poset whose probabilistic properties do not change under finite isomorphisms. They are Markov measures corresponding to  \textit{random walks on the poset}: such a measure is determined by the transition probabilities of choosing an element of~$P$ to be added to an ideal, which depend only on this ideal and the element to be added. The set of all ergodic central measures  will be denoted by
$\operatorname{Abs}(P)$ (absolute), it contains important information on the poset. From the viewpoint of measure theory, the problem of describing the list of central measures for posets has its own specific properties, which resemble the theory of invariant measures on groups.

An important characteristic of an ergodic central measure~$\mu$ is its
\textit{frequency function} on the space of ideals: $\Lambda_{\mu}(I)\equiv \lim_n
n^{-1}|\{i<n\colon \phi(i)\in I\}|$ for
$\mu$-a.e.\ numbering~$\phi$ and an infinite ideal~$I$. For a given measure~$\mu$, the frequency function~$\Lambda_\mu$ is a nonnegative monotone function on the space of infinite ideals which equals to~$1$ on the nonproper ideal coinciding with the entire poset and takes the same value on almost all numberings of every given central measure.

\begin{problem}
\label{prb1} Can an ergodic central measure be uniquely recovered from its frequency function? In particular, how many ergodic central measures are there such that the frequency function vanishes on all ideals except the nonproper nonzero ideal?
\end{problem}

For $P={\mathbb Z}_+^2$, an affirmative answer follows from a deep theorem of Thoma (see, e.g.,~\cite{2}) and, in particular, from its main corollary: there is a unique central measure (which was called the Plancherel measure) for which all proper infinite ideals have zero frequencies.

For $P={\mathbb Z}_+^d$, $d>2$, the problem is open.

The space of all infinite ideals of a poset is stratified by the dimension. The simplest ideals are one-dimensional ones: \textit{an ideal~$I$ is said to be one-dimensional if it has only finitely many incomparable elements}, in other words, if it is contained in the union of finitely many chains. \textit{An ergodic central measure~$\mu$ is said to be one-dimensional if its frequency function does not vanish only on one-dimensional ideals}.

The problem stated above is solved in the affirmative for one-dimensional central measures on the posets
$P={\mathbb Z}_+^d$, $d>1$.

\begin{theorem}
\label{th1} There is at most one ergodic central measure on the poset
$P={\mathbb Z}_+^d$, $d>1$, that has a given frequency function on the set of all one-dimensional ideals.
\end{theorem}

Already for quadratic
measures, proving a counterpart of the above theorem will require much more comprehensive information on the structure of Plancherel Young diagrams.

\end{document}